\documentclass[10pt]{article}
\begin{document}
\title{ {\bf  On Potentially $K_5-E_3$-graphic
Sequences}
\thanks{   Project Supported by NNSF of China(10271105), NSF of Fujian(Z0511034),
Science and Technology Project of Fujian, Fujian Provincial Training
Foundation for "Bai-Quan-Wan Talents Engineering" , Project of
Fujian Education Department and Project of Zhangzhou Teachers
College.}}
\author{{ Lili Hu , Chunhui Lai}\\
{\small Department of Mathematics, Zhangzhou Teachers College,}
\\{\small Zhangzhou, Fujian 363000,
 P. R. of CHINA.}\\{\small  jackey2591924@163.com ( Lili Hu)}
 \\{\small   zjlaichu@public.zzptt.fj.cn(Chunhui
 Lai, Corresponding author)}
}

\date{}
\maketitle
\begin{center}
\begin{minipage}{4.1in}
\vskip 0.1in
\begin{center}{\bf Abstract}\end{center}
 { Let $K_m-H$ be the graph obtained from $K_m$ by removing the edges
set $E(H)$ of $H$ where $H$ is a subgraph of $K_m$. In this paper,
we characterize the potentially $K_5-P_3$, $K_5-A_3$, $K_5-K_3$ and
$K_5-K_{1,3}$-graphic sequences where $A_3$ is $P_2\cup K_2$.
Moreover, we also characterize the potentially $K_5-2K_2$-graphic
sequences where $pK_2$ is the matching consisted of $p$ edges.}
\par
\end{minipage}
\end{center}
 \par
 \section{Introduction}
\par
\baselineskip 14pt
    We consider finite simple graphs. Any undefined notation follows
that of Bondy and Murty $[1]$. The set of all non-increasing
nonnegative integer sequence $\pi=(d_1,d_2,\cdots,d_n)$ is denoted
by $NS_n$. A sequence $\pi\epsilon NS_n$ is said to be graphic if it
is the degree sequence of a simple graph $G$ of order $n$; such a
graph $G$ is referred as a realization of $\pi$. The set of all
graphic sequence in $NS_n$ is denoted by $GS_n$. A graphic sequence
$\pi$ is potentially $H$-graphic if there is a realization of $\pi$
containing $H$ as a subgraph. Let $C_k$ and $P_k$ denote a cycle on
$k$ vertices and a path on $k+1$ vertices, respectively. Let
$\sigma(\pi)$ the sum of all the terms of $\pi$ and let $A_3$ and
$Z_4$ denote $P_2\cup K_2$ and $K_4-P_2$, respectively. We use the
symbol $E_3$ to denote graphs on 5 vertices and 3 edges. A graphic
sequence $\pi$ is said to be potentially $H$-graphic if it has a
realization $G$ containing $H$ as a subgraph. Let $G-H$ denote the
graph obtained from $G$ by removing the edges set $E(H)$ where $H$
is a subgraph of $G$.  In the degree sequence, $r^t$ means $r$
repeats $t$ times, that is, in the realization of the sequence there
are $t$ vertices of degree $r$.
\par

  Given a graph $H$, what is the maximum number of edges of a graph
with $n$ vertices not containing $H$ as a subgraph? This number is
denoted $ex(n,H)$, and is known as the Tur\'{a}n number. In terms of
graphic sequences, the number $2ex(n,H)+2$ is the minimum even
integer $l$ such that every $n$-term graphical sequence $\pi$ with
$\sigma (\pi)\geq l $ is forcibly $H$-graphical.  Gould, Jacobson
and Lehel [4] considered the following variation of the classical
Tur\'{a}n-type extremal problems: determine the smallest even
integer $\sigma(H,n)$ such that every n-term positive graphic
sequence $\pi=(d_1,d_2,\cdots,d_n)$ with $\sigma(\pi)\geq
\sigma(H,n)$ has a realization $G$ containing $H$ as a subgraph.
They proved that $\sigma(pK_2, n)=(p-1)(2n-p)+2$ for $p\ge 2$;
$\sigma(C_4, n)=2[{{3n-1}\over 2}]$ for $n\ge 4$.   Erd\"os,\
Jacobson and Lehel [3] showed that $\sigma(K_k, n)\ge
(k-2)(2n-k+1)+2$ and conjectured that the equality holds. In the
same paper, they proved the conjecture is true for $k=3$ and
$n\geq6$. The cases $k=4$ and 5 were proved separately (see [4] and
[17], and [18]). For $k\geq6$ and $n\geq {k \choose 2} $+3, Li, Song
and Luo [19] proved the conjecture true via linear algebraic
techniques. Recently, Ferrara, Gould and Schmitt  proved  the
conjecture $[5]$ and they also determined in $[6]$ $\sigma(F_k,n)$
where $F_k$ denotes the graph of $k$ triangles intersecting at
exactly one common vertex.
 Yin, Li, and Mao [25] determined $\sigma(K_{r+1}-e,n)$ for
$r\geq3$ and $r+1\leq n\leq 2r$ and $\sigma(K_5-e,n)$ for $n\geq5$,
and Yin and Li [24] further determined $\sigma(K_{r+1}-e,n)$ for
$r\geq2$ and $n\geq3r^2-r-1$. Moreover, Yin and Li in [24] also gave
two sufficient conditions for a sequence $\pi\epsilon GS_n$ to be
potentially $K_{r+1}-e$-graphic. Yin [27] determined
$\sigma(K_{r+1}-K_3,n)$ for $r\geq3$ and $n\geq 3r+5$. Lai [12-15]
determined  $\sigma(K_4-e,n)$ for $n\geq4$ and $\sigma(K_5-C_4,n)$,
$\sigma(K_5-P_3,n)$, $\sigma(K_5-P_4,n)$, $\sigma(K_5-K_3,n)$ for
$n\geq5$. Lai [10-11] proved that $\sigma(C_{2m+1}, n)=m(2n-m-1)+2$,
for $m\geq 2, n\geq 3m$; $\sigma(C_{2m+2} , n)=m(2n-m-1)+4$,
 for $ m\geq 2, n\geq 5m-2$.  Lai and Hu [16] determined $\sigma(K_{r+1}-H,n)$ for
$n\geq4r+10$, $r\geq3$, $r+1\geq k\geq4$ and $H$ be a graph on $k$
vertices which containing a tree on 4 vertices but not contain a
cycle on 3 vertices and $\sigma(K_{r+1}-P_2,n)$ for $n\geq4r+8$,
$r\geq3$.
\par
  A harder question is to characterize the potentially
 $H$-graphic sequences without zero terms.  Luo [21] characterized the potentially
 $C_k$-graphic sequences for each $k=3,4,5$. Recently, Luo and Warner [22] characterized the potentially
 $K_4$-graphic sequences.  Eschen and Niu [23] characterized the potentially
 $K_4-e$-graphic sequences.  Yin and Chen [26] characterized the
 potentially $K_{r,s}$-graphic sequences for $r=2,s=3$ and
 $r=2,s=4$. Chen [2] characterized the
 potentially $K_5-2K_2$-graphic sequences for $5\leq n\leq8$. Hu and Lai [7-8] characterized the potentially
 $K_5-C_4$ and $K_5-Z_4$-graphic
 sequences.
\par
  In this paper, we completely characterize the potentially $K_5-E_3$ -
graphic sequences, that is potentially $K_5-P_3$, $K_5-A_3$,
$K_5-K_3$ and $K_5-K_{1,3}$-graphic sequences. Moreover, we also
characterize the potentially $K_5-2K_2$-graphic sequences.
\par
\section{Preparations}\par
   Let $\pi=(d_1,\cdots,d_n)\epsilon NS_n,1\leq k\leq n$. Let
    $$ \pi_k^{\prime\prime}=\left\{
    \begin{array}{ll}(d_1-1,\cdots,d_{k-1}-1,d_{k+1}-1,
    \cdots,d_{d_k+1}-1,d_{d_k+2},\cdots,d_n), \\ \mbox{ if $d_k\geq k,$}\\
    (d_1-1,\cdots,d_{d_k}-1,d_{d_k+1},\cdots,d_{k-1},d_{k+1},\cdots,d_n),
     \\ \mbox{if $d_k < k.$} \end{array} \right. $$
  Denote
  $\pi_k^\prime=(d_1^\prime,d_2^\prime,\cdots,d_{n-1}^\prime)$, where
  $d_1^\prime\geq d_2^\prime\geq\cdots\geq d_{n-1}^\prime$ is a
  rearrangement of the $n-1$ terms of $\pi_k^{\prime\prime}$. Then
  $\pi_k^{\prime}$ is called the residual sequence obtained by
  laying off $d_k$ from $\pi$. In this paper, we denote $\pi_n^\prime$ by $\pi^\prime$.
\par
   For a nonincreasing positive integer sequence $\pi=(d_1,d_2,\cdots,d_n)$, we write $m(\pi)$ and $h(\pi)$ to denote the largest
positive terms of $\pi$ and the smallest positive terms of $\pi$,
respectively. We need the following results.
\par
    {\bf Theorem 2.1 [4]} If $\pi=(d_1,d_2,\cdots,d_n)$ is a graphic
 sequence with a realization $G$ containing $H$ as a subgraph,
 then there exists a realization $G^\prime$ of $\pi$ containing $H$ as a
 subgraph so that the vertices of $H$ have the largest degrees of
 $\pi$.\par
 \par
    {\bf Theorem 2.2 [20]} If $\pi=(d_1,d_2,\cdots,d_n)$ is a
 sequence of nonnegative integers with $1\leq m(\pi)\leq2$,
 $h(\pi)=1$ and even $\sigma(\pi)$, then $\pi$ is graphic.
\par
   {\bf Theorem 2.3 [21]} Let $\pi=(d_1,d_2,\cdots,d_n)$ be a graphic
 sequence. Then $\pi$ is potentially $C_4$-graphic if and only if
 the following conditions hold:
   (1) $d_4\geq2$;
   (2) $d_1=n-1$ implies $d_2\geq3$;
   (3) If $n=5,6$, then $\pi\neq(2^n)$.
\par
    {\bf Lemma 2.4 [2]} Let $\pi=(d_1,d_2,\cdots,d_n)\epsilon NS_n$,
$1\leq j\leq n-5$, $0\leq k\leq [{{n-j-i-4}\over 2}]$. Let $$ \pi =
\left\{
    \begin{array}{ll}(n-i,n-j,3^{n-i-j-2k}, 2^{2k},1^{i+j-2}) \\ \mbox{ $n-i-j$ is even;}\\
    (n-i,n-j,3^{n-i-j-2k-1}, 2^{2k+1},1^{i+j-2})
     \\ \mbox{ $n-i-j$ is odd.} \end{array} \right. $$
Let $S_1$ be the set consisting of the above sequences and let $S_2$
be the set of
  the following sequences: $(n-1,3^5,1^{n-6})$ and
  $(n-1,3^6,1^{n-7})$. If $\pi \epsilon S_1$ or $\pi \epsilon S_2$,
  then $\pi$ is not potentially $K_{1,2,2}$-graphic.
\par
    {\bf Lemma 2.5 [8]}  If $\pi=(d_1,d_2,\cdots,d_n)$ is a
 nonincreasing sequence of positive integers with even $\sigma(\pi)$,  $n\geq4$,
$d_1\leq3$ and $\pi\neq(3^3,1),(3^2,1^2)$, then $\pi$ is graphic.

 \par
    {\bf Lemma 2.6 (Kleitman and Wang [9])}\ \   $\pi$ is
graphic if and only if $\pi^\prime$ is graphic.
 \par
    The following corollary is obvious.\par
\par
    {\bf Corollary 2.7}\ \    Let $H$ be a simple graph. If $\pi^\prime$ is
 potentially $H$-graphic, then $\pi$ is
 potentially $H$-graphic.

\par
\section{ Main Theorems} \par
\par
\textbf{\noindent Theorem 3.1}  Let $\pi=(d_1,d_2,\cdots,d_n)$
 be a graphic sequence with $n\geq5$. Then $\pi$ is potentially
$K_5-P_3$-graphic if and only if the following conditions hold:
\par
  $(1)$ $d_1\geq4$, $d_3\geq3$ and $d_5\geq2$.
\par
  $(2)$ $\pi\neq (4,3^2,2^3)$, $(4,3^2,2^4)$ and $(4,3^6)$.
\par
{\bf Proof:} Assume that $\pi$ is potentially $K_5-P_3$-graphic.
$(1)$ and $(2)$ are obvious. To prove the sufficiency, we use
induction on $n$. Suppose the graphic sequence $\pi$ satisfies the
conditions (1) and (2). We first prove the base case where $n=5$. In
this case, $\pi$ is one of the following: $(4^5)$, $(4^3,3^2)$,
$(4^2,3^2,2)$, $(4,3^4)$, $(4,3^2,2^2)$. It is easy to check that
all of these are potentially $K_5-P_3$-graphic. Now we assume that
the sufficiency holds for $n-1(n\geq6)$, we will show that $\pi$ is
potentially $K_5-P_3$-graphic in terms of the following cases:
\par
\textbf{Case 1:} $d_n\geq4$. Clearly, $\pi^\prime$ satisfies $(1)$
and $(2)$, then by the induction hypothesis, $\pi^\prime$ is
potentially $K_5-P_3$-graphic, and hence so is $\pi$.
\par
\textbf{Case 2:} $d_n=3$. Consider
$\pi^\prime=(d_1^\prime,d_2^\prime,\cdots,d_{n-1}^\prime)$ where
$d_{n-3}^\prime\geq3$ and $d_{n-1}^\prime\geq2$. If $\pi^\prime$
satisfies $(1)$ and $(2)$, then by the induction hypothesis,
$\pi^\prime$ is potentially $K_5-P_3$-graphic, and hence so is
$\pi$.
\par
  If $\pi^\prime$ does not satisfy $(1)$, i.e., $d_1^\prime=3$, then $\pi^\prime=(3^k,2^{n-1-k})$ where
$n-3\leq k \leq n-1$. Since $\sigma(\pi^\prime)$ is even, $k$ must
be even.  If $k=n-3$, then $\pi=(4,3^{n-1})$ where $n$ is odd. Since
$\pi\neq(4,3^6)$, we have $n\geq9$. By Lemma 2.5, $\pi_1=(3^{n-5})$
is graphic. Let $G_1$ be a realization of $\pi_1$, then
$K_{1,2,2}\cup G_1$ is a realization of $\pi=(4,3^{n-1})$. Thus,
$\pi=(4,3^{n-1})$ is potentially $K_5-P_3$-graphic since
$K_5-P_3\subseteq K_{1,2,2}$. If $k=n-2$, then $\pi=(4^2,3^{n-2})$
where $n$ is even. It is easy to see that $\pi=(4^2,3^4)$ and
$\pi=(4^2,3^6)$ are potentially $K_5-P_3$-graphic. Let $G_2$ be a
realization of $(4^2,3^4)$, which contains $K_5-P_3$. If $n\geq10$,
then $\pi_2=(3^{n-6})$ is graphic by Lemma 2.5. Let $G_3$ be a
realization of $\pi_2$, then $G_2\cup G_3$ is a realization of
$\pi=(4^2,3^{n-2})$. In other words, $\pi=(4^2,3^{n-2})$ is
potentially $K_5-P_3$-graphic. If $k=n-1$, then $\pi=(4^3,3^{n-3})$
where $n$ is odd. It is easy to see that $\pi=(4^3,3^4)$ is
potentially $K_5-P_3$-graphic. If $n\geq9$, then $K_5-e\cup G_1$ is
a realization of $\pi=(4^3,3^{n-3})$. Thus, $\pi=(4^3,3^{n-3})$ is
potentially $K_5-P_3$-graphic since $K_5-P_3\subseteq K_5-e$.
\par
  If $\pi^\prime$ does not satisfy $(2)$, then $\pi^\prime$ is
just $(4,3^6)$, and hence $\pi=(5,4^2,3^5)$ or $(4^4,3^4)$. It is
easy to see that these sequences are potentially $K_5-P_3$-graphic.
\par
\textbf{Case 3:} $d_n=2$. Consider
$\pi^\prime=(d_1^\prime,d_2^\prime,\cdots,d_{n-1}^\prime)$ where
$d_2^\prime\geq3$ and $d_{n-1}^\prime\geq2$. If $\pi^\prime$
satisfies $(1)$ and $(2)$, then by the induction hypothesis,
$\pi^\prime$ is potentially $K_5-P_3$-graphic, and hence so is
$\pi$.
\par
  If $\pi^\prime$ does not satisfy $(1)$, there are two subcases:
\par
\textbf{Subcase 1:} $d_1^\prime\geq4$ and $d_3^\prime=2$. Then
$\pi=(d_1,3^2,2^{n-3})$ where $d_1\geq 5$. Since $\sigma(\pi)$ is
even, $d_1$ must be even. We will show that $\pi$ is potentially
$K_5-P_3$-graphic. It is enough to show $\pi_1=(d_1-4,2^{n-5})$ is
graphic. It clearly suffices to show $\pi_2=(2^{n-1-d_1},1^{d_1-4})$
is graphic. By $\sigma(\pi_2)$ being even and Theorem 2.2, $\pi_2$
is graphic.
\par
\textbf{Subcase 2:} $d_1^\prime=3$. Then $d_1=4$, $d_3=3$, $d_2=4$
or $d_2=3$.
\par
  If $d_2=4$, then $\pi=(4^2,3^k,2^{n-2-k})$ where $k\geq1$ and
$n-2-k\geq1$. Since $\sigma(\pi)$ is even, $k$ must be even. We will
show that $\pi$ is potentially $K_5-P_3$-graphic. First, we consider
$\pi=(4^2,3^2,2^{n-4})$. It is enough to show $\pi_1=(2^{n-5},1^2)$
is graphic.  By $\sigma(\pi_1)$ being even and Theorem 2.2, $\pi_1$
is graphic. Then we consider $\pi=(4^2,3^k,2^{n-2-k})$ where
$k\geq4$. It is easy to see that $(4^2,3^4)$ is potentially
$K_5-P_3$-graphic. Let $G_1$ be a realization of $(4^2,3^4)$, which
contains $K_5-P_3$. If $n\geq10$, then $\pi_2=(3^{k-4},2^{n-2-k})$
is graphic by Lemma 2.5. Let $G_2$ be a realization of $\pi_2$, then
$G_1\cup G_2$ is a realization of $\pi=(4^2,3^k,2^{n-2-k})$. If
$n\leq9$, then $\pi$ is one of the following: $(4^2,3^4,2)$,
$(4^2,3^4,2^2)$, $(4^2,3^4,2^3)$, $(4^2,3^6,2)$. It is easy to check
that all of these are potentially $K_5-P_3$-graphic. In other words,
$\pi=(4^2,3^k,2^{n-2-k})$ is potentially $K_5-P_3$-graphic.
\par
  If $d_2=3$, then $\pi=(4,3^k,2^{n-1-k})$ where $k\geq2$ and
$n-1-k\geq1$. Since $\sigma(\pi)$ is even, $k$ must be even. We will
show that $\pi$ is potentially $K_5-P_3$-graphic. First, we consider
$\pi=(4,3^2,2^{n-3})$. Since $\pi\neq(4,3^2,2^3)$ and $(4,3^2,2^4)$,
we have $n\geq8$.  It is enough to show $\pi_1=(2^{n-5})$ is
graphic. Clearly, $C_{n-5}$ is a realization of $\pi_1$. Second, we
consider $\pi=(4,3^4,2^{n-5})$. It is enough to show
$\pi_2=(2^{n-5},1^2)$ is graphic. By $\sigma(\pi_2)$ being even and
Theorem 2.2, $\pi_2$ is graphic. Then we consider
$\pi=(4,3^k,2^{n-1-k})$ where $k\geq6$. If $n\geq9$, then
$\pi_3=(3^{k-4},2^{n-1-k})$ is graphic by Lemma 2.5. Let $G_1$ be a
realization of $\pi_3$, then $K_{1,2,2}\cup G_1$ is a realization of
$\pi=(4,3^k,2^{n-1-k})$. Hence, $\pi=(4,3^k,2^{n-1-k})$ is
potentially $K_5-P_3$-graphic since $K_5-P_3\subseteq K_{1,2,2}$. If
$n\leq8$, then $\pi=(4,3^6,2)$. It is easy to see that $\pi$ is
potentially $K_5-P_3$-graphic. In other words,
$\pi=(4,3^k,2^{n-1-k})$ is potentially $K_5-P_3$-graphic.
\par
  If $\pi^\prime$ does not satisfy $(2)$, then $\pi^\prime$ is one
of the following: $(4,3^2,2^3)$, $(4,3^2,2^4)$, $(4,3^6)$. Hence
$\pi$ is one of the following: $(5,4,3,2^4)$, $(5,3^3,2^3)$,
$(4^3,2^4)$, $(5,4,3,2^5)$, $(5,3^3,2^4)$,  $(4^3,2^5)$,
$(5,4,3^5,2)$, $(4^3,3^4,2)$. It is easy to check that all of these
are potentially $K_5-P_3$-graphic.
\par
\textbf{Case 4:} $d_n=1$. Consider
$\pi^\prime=(d_1^\prime,d_2^\prime,\cdots,d_{n-1}^\prime)$ where
$d_3^\prime\geq3$ and $d_5^\prime\geq2$. If $\pi^\prime$ satisfies
$(1)$ and $(2)$, then by the induction hypothesis, $\pi^\prime$ is
potentially $K_5-P_3$-graphic, and hence so is $\pi$.
\par
  If $\pi^\prime$ does not satisfy $(1)$, i.e., $d_1^\prime=3$, then $\pi=(4,3^k,2^t,1^{n-1-k-t})$ where $k\geq2$, $k+t\geq4$ and
$n-1-k-t\geq1$. Since $\sigma(\pi)$ is even, $n-1-t$ must be even.
We will show that $\pi$ is potentially $K_5-P_3$-graphic. First, we
consider $\pi=(4,3^2,2^t,1^{n-3-t})$. It is enough to show
$\pi_1=(2^{t-2},1^{n-3-t})$ is graphic. By $\sigma(\pi_1)$ being
even and Theorem 2.2, $\pi_1$ is graphic.  Second, we consider
$\pi=(4,3^3,2^t,1^{n-4-t})$. It is enough to show
$\pi_2=(2^{t-1},1^{n-3-t})$ is graphic. By $\sigma(\pi_2)$ being
even and Theorem 2.2, $\pi_2$ is graphic. Third, we consider
$\pi=(4,3^4,2^t,1^{n-5-t})$. It is enough to show
$\pi_3=(2^t,1^{n-3-t})$ is graphic. By $\sigma(\pi_3)$ being even
and Theorem 2.2, $\pi_3$ is graphic. Then we consider
$\pi=(4,3^k,2^t,1^{n-1-k-t})$ where $k\geq5$. Let
$\pi_4=(3^{k-4},2^t,1^{n-1-k-t})$. If $n\geq9$ and
$\pi_4\neq(3^3,1)$ or $(3^2,1^2)$, then $\pi_4$ is graphic by Lemma
2.5. Let $G_1$ be a realization of $\pi_4$, then $K_{1,2,2}\cup G_1$
is a realization of $\pi=(4,3^k,2^t,1^{n-1-k-t})$. Hence,
$\pi=(4,3^k,2^t,1^{n-1-k-t})$ is potentially $K_5-P_3$-graphic since
$K_5-P_3\subseteq K_{1,2,2}$. If $n=9$ and $\pi_4=(3^3,1)$ or
$(3^2,1^2)$, then $\pi=(4,3^7,1)$ or $(4,3^6,1^2)$. If $n\leq8$,
then $\pi=(4,3^5,1)$ or $(4,3^5,2,1)$. It is easy to check that all
of these are potentially $K_5-P_3$-graphic. In other words,
$\pi=(4,3^k,2^t,1^{n-1-k-t})$ is potentially $K_5-P_3$-graphic.
\par
  If $\pi^\prime$ does not satisfy $(2)$, then $\pi^\prime$ is one
of the following: $(4,3^2,2^3)$, $(4,3^2,2^4)$, $(4,3^6)$. Hence
$\pi$ is one of the following: $(5,3^2,2^3,1)$, $(4^2,3,2^3,1)$,
$(5,3^2,2^4,1)$, $(4^2,3,2^4,1)$, $(5,3^6,1)$, $(4^2,3^5,1)$. It is
easy to check that all of these are potentially $K_5-P_3$-graphic.
\par
\vspace{0.5cm}
\par
\textbf{\noindent Theorem 3.2}  Let $\pi=(d_1,d_2,\cdots,d_n)$ be a
graphic sequence with $n\geq5$. Then $\pi$ is potentially
$K_5-A_3$-graphic if and only if the following conditions hold:
\par
  $(1)$ $d_4\geq3$ and $d_5\geq2$.
\par
  $(2)$ $\pi\neq(n-1,3^3,2^{n-k},1^{k-4})$ where $n\geq6$ and
$k=4,5,\cdots,n-2$, $n$ and $k$ have the same parity.
\par
  $(3)$ $\pi\neq (3^4,2^2),(3^6),(3^4,2^3),(3^6,2),(4,3^6),(3^7,1),(3^8),(n-1,3^5,1^{n-6})$
and $(n-1,3^6,1^{n-7})$.

\par
{\bf Proof:} First we show the conditions (1)-(3) are necessary
conditions for $\pi$ to be potentially $K_5-A_3$-graphic.  Assume
that $\pi$ is potentially $K_5-A_3$-graphic. $(1)$ is obvious. If
$\pi=(n-1,3^3,2^{n-k},1^{k-4})$ is potentially $K_5-A_3$-graphic,
then according to Theorem 2.1, there exists a realization $G$ of
$\pi$ containing $K_5-A_3$ as a subgraph so that the vertices of
$K_5-A_3$ have the largest degrees of $\pi$. Therefore, the sequence
$\pi^*=(n-4,2^{n-1-k},1^{k-4})$ obtained from $G-(K_5-A_3)$ must be
graphic, which is impossible since $G-(K_5-A_3)$ has only $n-4$
vertices, $\triangle(G-(K_5-A_3))\leq n-5$. Hence, $(2)$ holds. Now
it is easy to check that
$(3^4,2^2),(3^6),(3^4,2^3),(3^6,2),(4,3^6),(3^7,1)$ and $(3^8)$ are
not potentially $K_5-A_3$-graphic. If $\pi=(n-1,3^5,1^{n-6})$ is
potentially $K_5-A_3$-graphic, then according to Theorem 2.1, there
exists a realization $G$ of $\pi$ containing $K_5-A_3$ as a subgraph
so that the vertices of $K_5-A_3$ have the largest degrees of $\pi$.
Therefore, the sequence $\pi^*=(n-4,3,1^{n-5})$ obtained from
$G-(K_5-A_3)$ must be graphic. It follows that the sequence
$\pi_1=(2)$ must be graphic, a contradiction. Hence,
$\pi\neq(n-1,3^5,1^{n-6})$. If $\pi=(n-1,3^6,1^{n-7})$ is
potentially $K_5-A_3$-graphic, then according to Theorem 2.1, there
exists a realization $G$ of $\pi$ containing $K_5-A_3$ as a subgraph
so that the vertices of $K_5-A_3$ have the largest degrees of $\pi$.
Therefore, the sequence $\pi^*=(n-4,3^2,1^{n-6})$ obtained from
$G-(K_5-A_3)$ must be graphic. It follows that the sequence
$\pi_2=(2^2)$ must be graphic, a contradiction. Hence,
$\pi\neq(n-1,3^6,1^{n-7})$. In other words, $(3)$ holds.
\par
  Now we turn to show the conditions (1)-(3) are sufficient conditions for $\pi$
 to be potentially $K_5-A_3$-graphic. Suppose the graphic
sequence $\pi$ satisfies the conditions (1)-(3). Our proof is by
induction on $n$. We first prove the base case where $n=5$. In this
case, $\pi$ is one of the following: $(4^5)$, $(4^3,3^2)$,
$(4^2,3^2,2)$, $(4,3^4)$, $(3^4,2)$. It is easy to check that all of
these are potentially $K_5-A_3$-graphic. Now suppose that the
sufficiency holds for $n-1(n\geq6)$, we will show that $\pi$ is
potentially $K_5-A_3$-graphic in terms of the following cases:
\par
\textbf{Case 1:} $d_n\geq3$. Clearly, $\pi^\prime$ satisfies $(1)$.
If $\pi^\prime$ also satisfies $(2)$ and $(3)$, then by the
induction hypothesis, $\pi^\prime$ is potentially $K_5-A_3$-graphic,
and hence so is $\pi$.
\par
  If $\pi^\prime$ does not satisfy $(2)$, then $\pi^\prime$ is just
$(5,3^3,2^2)$, and hence $\pi=(6,3^6)$ which is impossible by $(3)$.
\par
  If $\pi^\prime$ does not satisfy $(3)$, since $\pi\neq(4,3^6)$ and $(3^8)$, then $\pi^\prime$ is
only one of the following: $(3^6),(3^6,2),(4,3^6)$, $(3^8)$,
$(5,3^5)$, $(6,3^6)$. Hence, $\pi$ is one of the following:
$(4^3,3^4),(4^2,3^6),(5,4^2,3^5),(4^4,3^4), (4^3,3^6)$,
$(6,4^2,3^4)$, $(7,4^2,3^5)$. It is easy to check that all of these
are potentially $K_5-A_3$-graphic.
\par
\textbf{Case 2:} $d_n=2$. Consider
$\pi^\prime=(d_1^\prime,d_2^\prime,\cdots,d_{n-1}^\prime)$ where
$d_2^\prime\geq3$ and $d_{n-1}^\prime\geq2$. If $\pi^\prime$
satisfies $(1)$-$(3)$, then by the induction hypothesis,
$\pi^\prime$ is potentially $K_5-A_3$-graphic, and hence so is
$\pi$.
\par
  If $\pi^\prime$ does not satisfy $(1)$, then $d_4^\prime=2$.
Hence $\pi=(d_1,3^3,2^{n-4})$. Since $\sigma(\pi)$ is even, $d_1$
must be odd. We will show that $\pi$ is potentially
$K_5-A_3$-graphic. If $d_1=3$, then $\pi=(3^4,2^{n-4})$. Since
$\pi\neq(3^4,2^2)$ and $(3^4,2^3)$, we have $n\geq8$. It is enough
to show $\pi_1=(2^{n-5})$ is graphic. Clearly, $C_{n-5}$ is a
realization of $\pi_1$. If $d_1\geq5$, since
$\pi\neq(n-1,3^3,2^{n-4})$, we have $d_1\leq n-2$. It is enough to
show $\pi_2=(d_1-3,2^{n-5})$ is graphic. It clearly suffices to show
$\pi_3=(2^{n-2-d_1},1^{d_1-3})$ is graphic. By $\sigma(\pi_3)$ being
even and Theorem 2.2, $\pi_3$ is graphic. Thus,
$\pi=(d_1,3^3,2^{n-4})$ is potentially $K_5-A_3$-graphic.
\par
  If $\pi^\prime$ does not satisfy $(2)$, i.e.,
$\pi^\prime=(n-2,3^3,2^{n-5})$.  Since $\sigma(\pi^\prime)$ is even,
$n$ must be odd. Hence $\pi=(n-1,4,3^2,2^{n-4})$ or
$(n-1,3^4,2^{n-5})$. We will show that both of them are potentially
$K_5-A_3$-graphic. It is enough to show $\pi_1=(n-4,2^{n-5},1)$ is
graphic. It clearly suffices to show $\pi_2=(1^{n-5})$ is graphic.
By $\sigma(\pi_2)$ being even and Theorem 2.2, $\pi_2$ is graphic.
\par
  If $\pi^\prime$ does not satisfy $(3)$, then $\pi^\prime$ is
one of the following: $(3^4,2^2)$, $(3^6)$, $(3^4,2^3)$, $(3^6,2)$,
$(4,3^6)$, $(3^8)$, $(5,3^5)$, $(6,3^6)$. Since $\pi\neq(3^6,2)$,
then $\pi$ is one of the following: $(4^2,3^2,2^3)$, $(4,3^4,2^2)$,
$(4^2,3^4,2)$, $(4^2,3^2,2^4)$, $(4,3^4,2^3)$, $(3^6,2^2)$,
$(4^2,3^4,2^2)$, $(4,3^6,2)$, $(5,4,3^5,2)$, $(4^3,3^4,2)$,
$(4^2,3^6,2)$, $(6,4,3^4,2)$, $(7,4,3^5,2)$. It is easy to check
that all of these are potentially $K_5-A_3$-graphic.
\par
\textbf{Case 3:} $d_n=1$. Consider
$\pi^\prime=(d_1^\prime,d_2^\prime,\cdots,d_{n-1}^\prime)$ where
$d_3^\prime\geq3$ and $d_5^\prime\geq2$. If $\pi^\prime$ satisfies
$(1)$-$(3)$, then by the induction hypothesis, $\pi^\prime$ is
potentially $K_5-A_3$-graphic, and hence so is $\pi$.
\par
  If $\pi^\prime$ does not satisfy $(1)$, then $d_4^\prime=2$.
Hence $\pi=(3^4,2^k,1^{n-4-k})$ where $k\geq1$ and $n-4-k\geq1$.
Since $\sigma(\pi)$ is even, $n-4-k$ must be even. We will show that
$\pi$ is potentially $K_5-A_3$-graphic. It is enough to show
$\pi_1=(2^{k-1},1^{n-4-k})$ is graphic. By $\sigma(\pi_1)$ being
even and Theorem 2.2, $\pi_1$ is graphic.
\par
  If $\pi^\prime$ does not satisfy $(2)$, i.e.,
$\pi^\prime=(n-2,3^3,2^{n-1-k},1^{k-4})$. Hence
$\pi=(n-1,3^3,2^{n-1-k},1^{k-3})$ which contradicts condition $(2)$.
\par
  If $\pi^\prime$ does not satisfy $(3)$, then by $\pi\neq(n-1,3^5,1^{n-6})$ and
$(n-1,3^6,1^{n-7})$, $\pi^\prime$ is only one of the following:
$(3^4,2^2)$, $(3^6)$, $(3^4,2^3)$, $(3^6,2)$, $(4,3^6)$, $(3^7,1)$,
$(3^8)$. Since $\pi\neq(3^7,1)$, then $\pi$ is one of the following:
$(4,3^3,2^2,1)$, $(3^5,2,1)$, $(4,3^5,1)$, $(4,3^3,2^3,1)$,
$(3^5,2^2,1)$, $(4,3^5,2,1)$, $(5,3^6,1)$,
 $(4^2,3^5,1)$, $(4,3^6,1^2)$, $(4,3^7,1)$.
 It is easy to check that all of these are potentially
$K_5-A_3$-graphic.
\par
 \vspace{0.5cm}
\par
\textbf{\noindent Theorem 3.3}  Let $\pi=(d_1,d_2,\cdots,d_n)$
 be a graphic sequence with $n\geq5$. Then $\pi$ is potentially
$K_5-K_3$-graphic if and only if the following conditions hold:
\par
  $(1)$ $d_2\geq4$ and $d_5\geq2$.
\par
  $(2)$ $\pi\neq (4^2,2^4)$, $(4^2,2^5)$, $(4^3,2^3)$ and $(4^6)$.
\par
{\bf Proof:} Assume that $\pi$ is potentially $K_5-K_3$-graphic.
$(1)$ and $(2)$ are obvious. To prove the sufficiency, we use
induction on $n$. Suppose the graphic sequence $\pi$ satisfies the
conditions (1) and (2). We first prove the base case where $n=5$. In
this case, $\pi$ is one of the following: $(4^5)$, $(4^3,3^2)$,
$(4^2,3^2,2)$, $(4^2,2^3)$. It is easy to check that all of these
are potentially $K_5-K_3$-graphic. Now suppose that the sufficiency
holds for $n-1(n\geq6)$, we will show that $\pi$ is potentially
$K_5-K_3$-graphic in terms of the following cases:
\par
\textbf{Case 1:} $d_n\geq4$. Clearly, $\pi^\prime$ satisfies $(1)$.
If $\pi^\prime$ also satisfies $(2)$, then by the induction
hypothesis, $\pi^\prime$ is potentially $K_5-K_3$-graphic, and hence
so is $\pi$. If $\pi^\prime$ does not satisfy $(2)$, then
$\pi^\prime$ is just $(4^6)$, and hence $\pi=(5^4,4^3)$. It is easy
to see that $\pi$ is potentially $K_5-K_3$-graphic.
\par
\textbf{Case 2:} $d_n=3$. Consider
$\pi^\prime=(d_1^\prime,d_2^\prime,\cdots,d_{n-1}^\prime)$ where
$d_{n-2}^\prime\geq3$ and $d_{n-1}^\prime\geq2$. If $\pi^\prime$
satisfies $(1)$ and $(2)$, then by the induction hypothesis,
$\pi^\prime$ is potentially $K_5-K_3$-graphic, and hence so is
$\pi$.
\par
  If $\pi^\prime$ does not satisfy $(1)$, i.e., $d_2^\prime=3$, then
$d_2=4$ and $3\leq d_4\leq d_3\leq4$. There are three subcases:
\par
\textbf{Subcase 1:} $d_4=4$. Then $\pi=(4^4,3^{n-4})$. Since
$\sigma(\pi)$ is even, $n$ must be even. We will show that $\pi$ is
potentially $K_5-K_3$-graphic. It is easy to see that $(4^4,3^2)$
and $(4^4,3^4)$ are potentially $K_5-K_3$-graphic. Let $G_1$ be a
realization of $(4^4,3^2)$, which contains $K_5-K_3$. If $n\geq10$,
then $\pi_1=(3^{n-6})$ is graphic by Lemma 2.5. Let $G_2$ be a
realization of $\pi_1$, then $G_1\cup G_2$ is a realization of
$\pi=(4^4,3^{n-4})$. In other words, $\pi=(4^4,3^{n-4})$ is
potentially $K_5-K_3$-graphic.

\par
\textbf{Subcase 2:} $d_4=3$ and $d_3=4$. Then
$\pi=(d_1,4^2,3^{n-3})$. Since $\sigma(\pi)$ is even, $d_1$ and $n$
have different parities. We will show that $\pi$ is potentially
$K_5-K_3$-graphic.  It is enough to show
$\pi_1=(d_1-4,3^{n-5},2,1^2)$ is graphic and the vertex with degree
$d_1-4$ is not adjacent to the vertices with degree 2 or 1 in the
realization of $\pi_1$. Hence, it suffices to show
$\pi_2=(3^{n-1-d_1},2^{d_1-3},1^2)$ is graphic. By Lemma 2.5,
$\pi_2$ is graphic. Thus, $\pi=(d_1,4^2,3^{n-3})$ is potentially
$K_5-K_3$-graphic.
\par
\textbf{Subcase 3:} $d_3=3$. then $\pi=(d_1,4,3^{n-2})$. Since
$\sigma(\pi)$ is even, $d_1$ and $n$ have the same parity. We will
show that $\pi$ is potentially $K_5-K_3$-graphic.  It is enough to
show $\pi_1=(d_1-4,3^{n-5},1^3)$ is graphic and the vertex with
degree $d_1-4$ is not adjacent to the vertices with degree 1 in the
realization of $\pi_1$. Hence, it suffices to show
$\pi_2=(3^{n-1-d_1},2^{d_1-4},1^3)$ is graphic. By Lemma 2.5,
$\pi_2$ is graphic.
\par
  If $\pi^\prime$ does not satisfy $(2)$, then $\pi^\prime$ is
just $(4^6)$, and hence $\pi=(5^3,4^3,3)$. It is easy to check that
$\pi$ is potentially $K_5-K_3$-graphic.
\par
\textbf{Case 3:} $d_n=2$. Consider
$\pi^\prime=(d_1^\prime,d_2^\prime,\cdots,d_{n-1}^\prime)$ where
$d_2^\prime\geq3$ and $d_{n-1}^\prime\geq2$. If $\pi^\prime$
satisfies $(1)$ and $(2)$, then by the induction hypothesis,
$\pi^\prime$ is potentially $K_5-K_3$-graphic, and hence so is
$\pi$.
\par
  If $\pi^\prime$ does not satisfy $(1)$, i.e., $d_2^\prime=3$, then $d_2=4$. There are two subcases:
$d_1=4$ and $d_1\geq5$.
\par
\textbf{Subcase 1:} $d_1=4$.
\par
  If $d_3=4$, then $\pi=(4^3,3^k,2^{n-3-k})$ where $n-3-k\geq1$. Since $\sigma(\pi)$
is even, $k$ must be even. We will show that $\pi$ is potentially
$K_5-K_3$-graphic. First, we consider $\pi=(4^3,2^{n-3})$. Since
$\pi\neq(4^3,2^3)$, we have $n\geq7$. It is enough to show
$\pi_1=(2^{n-4})$ is graphic. Clearly, $C_{n-4}$ is a realization of
$\pi_1$. Second, we consider $\pi=(4^3,3^2,2^{n-5})$.  It is easy to
see that $\pi=(4^3,3^2,2)$ and $\pi=(4^3,3^2,2^2)$ are potentially
$K_5-K_3$-graphic. If $n\geq8$, then $K_5-e \cup C_{n-5}$ is a
realization of $\pi=(4^3,3^2,2^{n-5})$. Thus,
$\pi=(4^3,3^2,2^{n-5})$ is potentially $K_5-K_3$-graphic since
$K_5-K_3\subseteq K_5-e$. Then we consider $\pi=(4^3,3^k,2^{n-3-k})$
where $k\geq4$. If $n\geq9$, then $\pi_1=(3^{k-2},2^{n-3-k})$ is
graphic by Lemma 2.5. Let $G_1$ be a realization of $\pi_1$, then
$K_5-e \cup G_1$ is a realization of $\pi=(4^3,3^k,2^{n-3-k})$.
Thus, $\pi$ is potentially $K_5-K_3$-graphic since $K_5-K_3\subseteq
K_5-e$. If $n\leq8$, then $\pi=(4^3,3^4,2)$. It is easy to see that
$(4^3,3^4,2)$ is potentially $K_5-K_3$-graphic. In other words,
$\pi=(4^3,3^k,2^{n-3-k})$ is potentially $K_5-K_3$-graphic.
\par
  If $d_3\leq3$, then $\pi=(4^2,3^k,2^{n-2-k})$ where $n-2-k\geq1$. Since $\sigma(\pi)$
is even, $k$ must be even. We will show that $\pi$ is potentially
$K_5-K_3$-graphic. First, we consider $\pi=(4^2,2^{n-2})$. Since
$\pi\neq(4^2,2^4)$ and $(4^2,2^5)$, we have $n\geq8$. It is enough
to show $\pi_1=(2^{n-5})$ is graphic. Clearly, $C_{n-5}$ is a
realization of $\pi_1$. Second, we consider $\pi=(4^2,3^2,2^{n-4})$.
It is enough to show $\pi_2=(2^{n-5},1^2)$ is graphic.  By
$\sigma(\pi_2)$ being even and Theorem 2.2, $\pi_2$ is graphic. Then
we consider $\pi=(4^2,3^k,2^{n-2-k})$ where $k\geq4$. It is easy to
check that $\pi=(4^2,3^4)$ is potentially $K_5-K_3$-graphic. Let
$G_1$ be a realization of $(4^2,3^4)$, which contains $K_5-K_3$. If
$n\geq10$, then $\pi_3=(3^{k-4},2^{n-2-k})$ is graphic by Lemma 2.5.
Let $G_2$ be a realization of $\pi_3$, then $G_1\cup G_2$ is a
realization of $\pi=(4^2,3^k,2^{n-2-k})$. If $n\leq9$, then $\pi$ is
one of the following: $(4^2,3^4,2)$, $(4^2,3^4,2^2)$,
$(4^2,3^4,2^3)$, $(4^2,3^6,2)$. It is easy to check that all of
these are potentially $K_5-K_3$-graphic. In other words,
$\pi=(4^2,3^k,2^{n-2-k})$ is potentially $K_5-K_3$-graphic.
\par
\textbf{Subcase 2:} $d_1\geq5$. Then $\pi=(d_1,4,3^k,2^{n-2-k})$
where $n-2-k\geq1$. Since $\sigma(\pi)$ is even, $d_1$ and $k$ have
the same parity. We will show that $\pi$ is potentially
$K_5-K_3$-graphic.
\par
  First, we consider $\pi=(d_1,4,2^{n-2})$.  It is enough
to show $\pi_1=(d_1-4,2^{n-5})$ is graphic. It clearly suffices to
show $\pi_2=(2^{n-1-d_1},1^{d_1-4})$ is graphic. By $\sigma(\pi_2)$
being even and Theorem 2.2, $\pi_2$ is graphic.
\par
  Second, we consider $\pi=(d_1,4,3,2^{n-3})$.  It is enough
to show $\pi_1=(d_1-4,2^{n-5},1)$ is graphic and there exists no
edge between two vertices with degree $d_1-4$ and $1$ in the
realization of $\pi_1$. Hence, it suffices to show
$\pi_2=(2^{n-1-d_1},1^{d_1-3})$ is graphic. By $\sigma(\pi_2)$ being
even and Theorem 2.2, $\pi_2$ is graphic.
\par
  Third, we consider $\pi=(d_1,4,3^2,2^{n-4})$.  It is enough
to show $\pi_1=(d_1-4,2^{n-5},1^2)$ is graphic and the vertex with
degree $d_1-4$ is not adjacent to the vertices with degree $1$ in
the realization of $\pi_1$. Hence, it suffices to show
$\pi_2=(2^{n-1-d_1},1^{d_1-2})$ is graphic. By $\sigma(\pi_2)$ being
even and Theorem 2.2, $\pi_2$ is graphic.
\par
  Fourth, we consider $\pi=(d_1,4,3^3,2^{n-5})$.  It is enough
to show $\pi_1=(d_1-4,2^{n-5},1^3)$ is graphic and the vertex with
degree $d_1-4$ is not adjacent to the vertices with degree $1$ in
the realization of $\pi_1$. Hence, it suffices to show
$\pi_2=(2^{n-1-d_1},1^{d_1-1})$ is graphic. By $\sigma(\pi_2)$ being
even and Theorem 2.2, $\pi_2$ is graphic.
\par
  Then we consider $\pi=(d_1,4,3^k,2^{n-2-k})$ where $k\geq4$.  It is enough
to show $\pi_1=(d_1-4,3^{k-3},2^{n-2-k},1^3)$ is graphic and the
vertex with degree $d_1-4$ is not adjacent to the vertices with
degree $1$ in the realization of $\pi_1$.  Assume that the vertex
with degree $d_1-4$ is adjacent to $t$$(t\leq k-3)$ vertices with
degree $3$ and $d_1-4-t$ vertices with degree $2$ in the realization
of $\pi_1$. Hence, it suffices to show
$\pi_2=(3^{k-3-t},2^{n+2-d_1-k+2t},1^{d_1-1-t})$ is graphic. By
Lemma 2.5, $\pi_2$ is graphic. Thus, $\pi=(d_1,4,3^k,2^{n-2-k})$ is
potentially $K_5-K_3$-graphic.
\par
  If $\pi^\prime$ does not satisfy $(2)$, then $\pi^\prime$ is
one of the following: $(4^2,2^4)$, $(4^2,2^5)$, $(4^3,2^3)$,
$(4^6)$. Hence $\pi$ is one of the following: $(5^2,2^5)$,
$(5^2,2^6)$, $(5^2,4,2^4)$, $(5^2,4^4,2)$. It is easy to check that
all of these are potentially $K_5-K_3$-graphic.
\par
\textbf{Case 4:} $d_n=1$. Consider
$\pi^\prime=(d_1^\prime,d_2^\prime,\cdots,d_{n-1}^\prime)$ where
$d_1^\prime\geq4$, $d_2^\prime\geq3$ and $d_5^\prime\geq2$. If
$\pi^\prime$ satisfies $(1)$ and $(2)$, then by the induction
hypothesis, $\pi^\prime$ is potentially $K_5-K_3$-graphic, and hence
so is $\pi$.
\par
  If $\pi^\prime$ does not satisfy $(1)$, i.e., $d_2^\prime=3$,
then $\pi=(4^2,3^k,2^t,1^{n-2-k-t})$ where $k+t\geq3$ and
$n-2-k-t\geq1$. Since $\sigma(\pi)$ is even, $n-2-t$ must be even.
We will show that $\pi$ is potentially $K_5-K_3$-graphic.
\par
  First, we consider $\pi=(4^2,2^t,1^{n-2-t})$.  It is enough
to show $\pi_1=(2^{t-3},1^{n-2-t})$ is graphic. By $\sigma(\pi_1)$
being even and Theorem 2.2, $\pi_1$ is graphic.
\par
  Second, we consider $\pi=(4^2,3,2^t,1^{n-3-t})$.  It is enough
to show $\pi_1=(2^{t-2},1^{n-2-t})$ is graphic. By $\sigma(\pi_1)$
being even and Theorem 2.2, $\pi_1$ is graphic.
\par
  Third, we consider $\pi=(4^2,3^2,2^t,1^{n-4-t})$.  It is enough
to show $\pi_1=(2^{t-1},1^{n-2-t})$ is graphic. By $\sigma(\pi_1)$
being even and Theorem 2.2, $\pi_1$ is graphic.
\par
  Fourth, we consider $\pi=(4^2,3^3,2^t,1^{n-5-t})$.  It is enough
to show $\pi_1=(2^t,1^{n-2-t})$ is graphic. By $\sigma(\pi_1)$ being
even and Theorem 2.2, $\pi_1$ is graphic.
\par
  Then we consider $\pi=(4^2,3^k,2^t,1^{n-2-k-t})$ where $k\geq4$ and
$n-2-k-t\geq1$. It is easy to see that $\pi=(4^2,3^4)$ is
potentially $K_5-K_3$-graphic. Let $G_1$ be a realization of
$(4^2,3^4)$, which contains $K_5-K_3$. Let
$\pi_1=(3^{k-4},2^t,1^{n-2-k-t})$. If $n\geq10$ and
$\pi_1\neq(3^3,1)$, $(3^2,1^2)$, then $\pi_1$ is graphic by Lemma
2.5. Let $G_2$ be a realization of $\pi_1$, then $G_1 \cup G_2$ is a
realization of $\pi=(4^2,3^k,2^t,1^{n-2-k-t})$. If $n=10$ and
$\pi_1=(3^3,1)$ or $(3^2,1^2)$, then $\pi=(4^2,3^7,1)$ or
$(4^2,3^6,1^2)$. If $n\leq9$, then $\pi=(4^2,3^4,1^2)$,
$(4^2,3^4,2,1^2)$, $(4^2,3^5,1)$ or $(4^2,3^5,2,1)$.  It is easy to
check that all of these are potentially $K_5-K_3$-graphic. In other
words, $\pi=(4^2,3^k,2^t,1^{n-2-k-t})$ is potentially
$K_5-K_3$-graphic.
\par
  If $\pi^\prime$ does not satisfy $(2)$, then $\pi^\prime$ is
one of the following: $(4^2,2^4)$, $(4^2,2^5)$, $(4^3,2^3)$,
$(4^6)$. Hence $\pi$ is one of the following: $(5,4,2^4,1)$,
$(5,4,2^5,1)$,  $(5,4^2,2^3,1)$, $(5,4^5,1)$. It is easy to check
that all of these are potentially $K_5-K_3$-graphic.
\par
\vspace{0.5cm}
\par
\textbf{\noindent Theorem 3.4}  Let $\pi=(d_1,d_2,\cdots,d_n)$ be a
graphic sequence with $n\geq5$. Then $\pi$ is potentially
$K_5-K_{1,3}$-graphic if and only if the following conditions hold:
\par
  $(1)$ $d_1\geq4$ and $d_4\geq3$.
\par
  $(2)$ $\pi\neq (4,3^4,2)$, $(4^6)$, $(4^2,3^4)$, $(4,3^6)$, $(4^7)$, $(4,3^5,1)$, $(n-1,3^4,1^{n-5})$
and $(n-1,3^5,1^{n-6})$.

\par
{\bf Proof:} Assume that $\pi$ is potentially $K_5-K_{1,3}$-graphic.
$(1)$ is obvious.  Now it is easy to check that $(4,3^4,2)$,
$(4^6)$, $(4^2,3^4)$, $(4,3^6)$, $(4^7)$, $(4,3^5,1)$ are not
potentially $K_5-K_{1,3}$-graphic. If $\pi=(n-1,3^4,1^{n-5})$ is
potentially $K_5-K_{1,3}$-graphic, then according to Theorem 2.1,
there exists a realization $G$ of $\pi$ containing $K_5-K_{1,3}$ as
a subgraph so that the vertices of $K_5-K_{1,3}$ have the largest
degrees of $\pi$. Therefore, the sequence $\pi^*=(n-5,2,1^{n-5})$
obtained from $G-(K_5-K_{1,3})$ must be graphic and there must be no
edge between two vertices with degree $n-5$ and $2$ in the
realization of $\pi^*$. Thus, $\pi^*$ satisfies: $(n-5)+2\leq n-5$,
a contradiction. Hence, $\pi\neq(n-1,3^4,1^{n-5})$.  If
$\pi=(n-1,3^5,1^{n-6})$ is potentially $K_5-K_{1,3}$-graphic, then
according to Theorem 2.1, there exists a realization $G$ of $\pi$
containing $K_5-K_{1,3}$ as a subgraph so that the vertices of
$K_5-K_{1,3}$ have the largest degrees of $\pi$. Therefore, the
sequence $\pi^*=(n-5,3,2,1^{n-6})$ obtained from $G-(K_5-K_{1,3})$
must be graphic and there must be no edge between two vertices with
degree $n-5$ and $2$ in the realization of $\pi^*$. It follows that
the sequence $\pi_1=(2^2)$ must be graphic, a contradiction. Hence,
$\pi\neq(n-1,3^5,1^{n-6})$. In other words, $(2)$ holds.
\par
  Now we prove the sufficient conditions. Suppose the graphic
sequence $\pi$ satisfies the conditions $(1)$ and $(2)$. Our proof
is by induction on $n$. We first prove the base case where $n=5$.
Since $\pi\neq(4,3^4)$, then $\pi$ is one of the following: $(4^5)$,
$(4^3,3^2)$, $(4^2,3^2,2)$, $(4,3^3,1)$. It is easy to check that
all of these are potentially $K_5-K_{1,3}$-graphic. Now suppose that
the sufficiency holds for $n-1(n\geq6)$, we will show that $\pi$ is
potentially $K_5-K_{1,3}$-graphic in terms of the following cases:
\par
\textbf{Case 1:} $d_n\geq4$. Clearly, $\pi^\prime$ satisfies $(1)$.
If $\pi^\prime$ also satisfies $(2)$, then by the induction
hypothesis, $\pi^\prime$ is potentially $K_5-K_{1,3}$-graphic, and
hence so is $\pi$. If $\pi^\prime$ does not satisfy $(2)$, since
$\pi\neq(4^6)$ and $(4^7)$, then $\pi^\prime$ is just $(4^6)$ or
$(4^7)$, and hence $\pi=(5^4,4^3)$ or $(5^4,4^4)$. It is easy to
check that these sequences are potentially $K_5-K_{1,3}$-graphic.
\par
\textbf{Case 2:} $d_n=3$. Consider
$\pi^\prime=(d_1^\prime,d_2^\prime,\cdots,d_{n-1}^\prime)$ where
$d_{n-3}^\prime\geq3$ and $d_{n-1}^\prime\geq2$. If $\pi^\prime$
satisfies $(1)$ and $(2)$, then by the induction hypothesis,
$\pi^\prime$ is potentially $K_5-K_{1,3}$-graphic, and hence so is
$\pi$.
\par
  If $\pi^\prime$ does not satisfy $(1)$, there are two subcases:
\par
\textbf{Subcase 1:} $d_1^\prime\geq4$ and $d_4^\prime=2$. Then
$\pi^\prime=(4,3^2,2^2)$, and hence $\pi=(5,3^5)$ which
contradicts condition $(2)$.
\par
\textbf{Subcase 2:} $d_1^\prime=3$. Then
$\pi^\prime=(3^k,2^{n-1-k})$ where $n-3\leq k\leq n-1$. Since
$\sigma(\pi^\prime)$ is even, $k$ must be even. If $n$ is odd, then
$k=n-3$ or $n-1$. If $k=n-3$, then $\pi=(4,3^{n-1})$. Since
$\pi\neq(4,3^6)$, we have $n\geq9$. It is easy to check that
$(4,3^8)$ and $(4,3^{10})$ are potentially $K_5-K_{1,3}$-graphic.
Let $G_1$ be a realization of $(4,3^8)$, which contains
$K_5-K_{1,3}$. If $n\geq13$, then $\pi_1=(3^{n-9})$ is graphic by
Lemma 2.5. Let $G_2$ be a realization of $\pi_1$, then $G_1\cup G_2$
is a realization of $\pi=(4,3^{n-1})$. In other words,
$\pi=(4,3^{n-1})$ is potentially $K_5-K_{1,3}$-graphic. If $k=n-1$,
then $\pi=(4^3,3^{n-3})$. It is easy to see that $\pi=(4^3,3^4)$ is
potentially $K_5-K_{1,3}$-graphic. If $n\geq9$, then
$\pi_2=(3^{n-5})$ is graphic by Lemma 2.5. Let $G_3$ be a
realization of $\pi_2$, then $K_5-e\cup G_3$ is a realization of
$\pi=(4^3,3^{n-3})$. Hence, $\pi=(4^3,3^{n-3})$ is potentially
$K_5-K_{1,3}$-graphic since $K_5-K_{1,3}\subseteq K_5-e$. If $n$ is
even, then $k=n-2$, thus $\pi=(4^2,3^{n-2})$. Since $\pi\neq
(4^2,3^4)$, we have $n\geq8$. It is easy to see that $(4^2,3^6)$ and
$(4^2,3^8)$ are potentially $K_5-K_{1,3}$-graphic. Let $G_4$ be a
realization of $(4^2,3^6)$, which contains $K_5-K_{1,3}$. If
$n\geq12$, then $\pi_3=(3^{n-8})$ is graphic by Lemma 2.5. Let $G_5$
be a realization of $\pi_3$, then $G_4\cup G_5$ is a realization of
$\pi=(4^2,3^{n-2})$. In other words, $\pi=(4^2,3^{n-2})$ is
potentially $K_5-K_{1,3}$-graphic.
\par
  If $\pi^\prime$ does not satisfy $(2)$, then $\pi^\prime$ is
one of the following: $(4,3^4,2)$, $(4^6)$, $(4^2,3^4)$, $(4,3^6)$,
$(4^7)$,  $(4,3^4)$, $(5,3^5)$. Hence $\pi$ is one of the following:
$(5,4,3^5)$, $(5^3,4^3,3)$,  $(5^2,4,3^4)$, $(5,4^3,3^3)$,
$(4^5,3^2)$, $(5,4^2,3^5)$, $(4^4,3^4)$, $(5^3,4^4,3)$,
$(5,4^2,3^3)$, $(4^4,3^2)$, $(6,4^2,3^4)$. It is easy to check that
all of these are potentially $K_5-K_{1,3}$-graphic.
\par
\textbf{Case 3:} $d_n=2$. Consider
$\pi^\prime=(d_1^\prime,d_2^\prime,\cdots,d_{n-1}^\prime)$ where
$d_3^\prime\geq3$ and $d_{n-1}^\prime\geq2$. If $\pi^\prime$
satisfies $(1)$ and $(2)$, then by the induction hypothesis,
$\pi^\prime$ is potentially $K_5-K_{1,3}$-graphic, and hence so is
$\pi$.
\par
  If $\pi^\prime$ does not satisfy $(1)$, there are two subcases:
\par
\textbf{Subcase 1:} $d_1^\prime\geq4$ and $d_4^\prime=2$. Then
$\pi=(d_1,3^3,2^{n-4})$ where $d_1\geq5$. Since $\sigma(\pi)$ is
even, $d_1$ must be odd. We will show that $\pi$ is potentially
$K_5-K_{1,3}$-graphic. It is enough to show
$\pi_1=(d_1-4,2^{n-5},1)$ is graphic and there exists no edge
between two vertices with degree $d_1-4$ and $1$ in the realization
of $\pi_1$. Hence, it suffices to show
$\pi_2=(2^{n-1-d_1},1^{d_1-3})$ is graphic. By $\sigma(\pi_2)$ being
even and Theorem 2.2, $\pi_2$ is graphic.
\par
\textbf{Subcase 2:} $d_1^\prime=3$. Then $d_1=4$, $d_3=d_4=3$,
$d_2=4$ or $d_2=3$.
\par
  If $d_2=4$, then $\pi=(4^2,3^k,2^{n-2-k})$ where $k\geq2$ and
$n-2-k\geq1$.  Since $\sigma(\pi)$ is even, $k$ must be even. We
will show that $\pi$ is potentially $K_5-K_{1,3}$-graphic. First, we
consider $\pi=(4^2,3^2,2^{n-4})$. It is enough to show
$\pi_1=(2^{n-5},1^2)$ is graphic. By $\sigma(\pi_1)$ being even and
Theorem 2.2, $\pi_1$ is graphic. Second, we consider
$\pi=(4^2,3^4,2^{n-6})$. It is easy to see that $(4^2,3^4,2)$,
$(4^2,3^4,2^2)$ and $(4^2,3^4,2^3)$ are potentially
$K_5-K_{1,3}$-graphic. Let $G_1$ be a realization of $(4^2,3^4,2)$,
which contains $K_5-K_{1,3}$. If $n\geq10$, then
 $G_1\cup C_{n-7}$ is a realization of
$\pi=(4^2,3^4,2^{n-6})$. In other words, $\pi=(4^2,3^4,2^{n-6})$ is
potentially $K_5-K_{1,3}$-graphic. Then we consider
$\pi=(4^2,3^k,2^{n-2-k})$ where $k\geq6$. It is easy to see that
$\pi=(4^2,3^6)$ is potentially $K_5-K_{1,3}$-graphic.  Let $G_2$ be
a realization of $(4^2,3^6)$, which contains $K_5-K_{1,3}$. If
$n\geq12$, then $\pi_2=(3^{k-6},2^{n-2-k})$ is graphic by Lemma 2.5.
Let $G_3$ be a realization of $\pi_2$, then $G_2\cup G_3$ is a
realization of $\pi=(4^2,3^k,2^{n-2-k})$. If $n\leq11$, then $\pi$
is one of the following: $(4^2,3^6,2)$, $(4^2,3^6,2^2)$,
$(4^2,3^6,2^3)$, $(4^2,3^8,2)$. It is easy to check that all of
these are potentially $K_5-K_{1,3}$-graphic. In other words,
$\pi=(4^2,3^k,2^{n-2-k})$ is potentially $K_5-K_{1,3}$-graphic.
\par
  If $d_2=3$, then $\pi=(4,3^k,2^{n-1-k})$ where $k\geq3$ and
$n-1-k\geq1$.  Since $\sigma(\pi)$ is even, $k$ must be even. We
will show that $\pi$ is potentially $K_5-K_{1,3}$-graphic. First, we
consider $\pi=(4,3^4,2^{n-5})$. Since $\pi\neq(4,3^4,2)$, we have
$n\geq7$. It is enough to show $\pi_1=(2^{n-4})$ is graphic.
Clearly, $C_{n-4}$ is a realization of $\pi_1$. Second, we consider
$\pi=(4,3^6,2^{n-7})$. It is easy to see that $(4,3^6,2)$,
$(4,3^6,2^2)$ and $(4,3^6,2^3)$ are potentially
$K_5-K_{1,3}$-graphic.  Let $G_1$ be a realization of $(4,3^6,2)$,
which contains $K_5-K_{1,3}$. If $n\geq11$, then
 $G_1\cup C_{n-8}$ is a realization of
$\pi=(4,3^6,2^{n-7})$. In other words, $\pi=(4,3^6,2^{n-7})$ is
potentially $K_5-K_{1,3}$-graphic. Then we consider
$\pi=(4,3^k,2^{n-1-k})$ where $k\geq8$. It is easy to see that
$\pi=(4^,3^8)$ is potentially $K_5-K_{1,3}$-graphic.  Let $G_2$ be a
realization of $(4,3^8)$, which contains $K_5-K_{1,3}$. If
$n\geq13$, then $\pi_2=(3^{k-8},2^{n-1-k})$ is graphic by Lemma 2.5.
Let $G_3$ be a realization of $\pi_2$, then $G_2\cup G_3$ is a
realization of $\pi=(4,3^k,2^{n-1-k})$. If $n\leq12$, then $\pi$ is
one of the following: $(4,3^8,2)$, $(4,3^8,2^2)$, $(4,3^8,2^3)$,
$(4,3^{10},2)$. It is easy to check that all of these are
potentially $K_5-K_{1,3}$-graphic. In other words,
$\pi=(4,3^k,2^{n-1-k})$ is potentially $K_5-K_{1,3}$-graphic.
\par
  If $\pi^\prime$ does not satisfy $(2)$, then $\pi^\prime$ is
one of the following: $(4,3^4,2)$, $(4^6)$, $(4^2,3^4)$,
  $(4,3^6)$, $(4^7)$, $(4,3^4)$, $(5,3^5)$. Hence $\pi$ is one of the following:
 $(5,4,3^3,2^2)$, $(5,3^5,2)$,
$(4^3,3^2,2^2)$,  $(5^2,4^4,2)$, $(5^2,3^4,2)$, $(5,4^2,3^3,2)$,
$(4^4,3^2,2)$, $(5,4,3^5,2)$, $(4^3,3^4,2)$, $(5^2,4^5,2)$,
$(5,4,3^3,2)$, $(4^3,3^2,2)$, $(6,4,3^4,2)$. It is easy to check
that all of these are potentially $K_5-K_{1,3}$-graphic.
\par
\textbf{Case 4:} $d_n=1$. Consider
$\pi^\prime=(d_1^\prime,d_2^\prime,\cdots,d_{n-1}^\prime)$ where
$d_4^\prime\geq3$. If $\pi^\prime$ satisfies $(1)$ and $(2)$, then
by the induction hypothesis, $\pi^\prime$ is potentially
$K_5-K_{1,3}$-graphic, and hence so is $\pi$.
\par
  If $\pi^\prime$ does not satisfy $(1)$, i.e.,
$d_1^\prime=3$, then $\pi=(4,3^k,2^t,1^{n-1-k-t})$ where $k\geq3$
and $n-1-k-t\geq1$. Since $\sigma(\pi)$ is even, $n-1-t$ must be
even. We will show that $\pi$ is potentially $K_5-K_{1,3}$-graphic.
\par
  First, we consider $\pi=(4,3^3,2^t,1^{n-4-t})$. If $t=0$, it is enough to show
$\pi_1=(1^{n-5})$ is graphic. By $\sigma(\pi_1)$ being even and
Theorem 2.2, $\pi_1$ is graphic. If $t\geq1$, it is enough to show
$\pi_2=(2^{t-1},1^{n-3-t})$ is graphic. By $\sigma(\pi_2)$ being
even and Theorem 2.2, $\pi_2$ is graphic.
\par
  Second, we consider $\pi=(4,3^4,2^t,1^{n-5-t})$. It is enough to show
$\pi_1=(2^{t+1},1^{n-5-t})$ is graphic. By $\sigma(\pi_1)$ being
even and Theorem 2.2, $\pi_1$ is graphic.
\par
  Then we consider $\pi=(4,3^k,2^t,1^{n-1-k-t})$ where $k\geq5$.
Since $\pi\neq(4,3^5,1)$, we have $n\geq8$. It is enough to show
$\pi_1=(3^{k-4},2^{t+1},1^{n-1-k-t})$ is graphic. By Lemma 2.5,
$\pi_1$ is graphic.
\par
  If $\pi^\prime$ does not satisfy $(2)$, since $\pi\neq(n-1,3^4,1^{n-5})$ and $(n-1,3^5,1^{n-6})$,
then $\pi^\prime$ is one of the following: $(4,3^4,2)$, $(4^6)$,
$(4^2,3^4)$,
 $(4,3^6)$, $(4^7)$, $(4,3^5,1)$. Hence,  $\pi$ is one of the
 following: $(5,3^4,2,1)$,
$(4^2,3^3,2,1)$, $(5,4^5,1)$, $(5,4,3^4,1)$,  $(4^3,3^3,1)$,
$(5,3^6,1)$, $(4^2,3^5,1)$, $(5,4^6,1)$, $(5,3^5,1^2)$,
$(4^2,3^4,1^2)$. It is easy to check that all of these are
potentially $K_5-K_{1,3}$-graphic.
\par
\vspace{0.5cm}
\par
\textbf{\noindent Theorem 3.5}  Let $\pi=(d_1,d_2,\cdots,d_n)$ be a
graphic sequence with $n\geq5$. Then $\pi$ is potentially
$K_5-2K_2$-graphic if and only if the following conditions hold:
\par
  $(1)$ $d_1\geq4$ and $d_5\geq3$;
\par
  $(2)$ $$ \pi \neq \left\{
    \begin{array}{ll}(n-i,n-j,3^{n-i-j-2k}, 2^{2k},1^{i+j-2})\\  \mbox{ $n-i-j$ is even;}\\
    (n-i,n-j,3^{n-i-j-2k-1}, 2^{2k+1},1^{i+j-2})
  \\     \mbox{$n-i-j$ is odd.} \end{array} \right. $$ \ \ \ \ \ where $1\leq j\leq
     n-5$ and $0\leq k\leq [{{n-j-i-4}\over 2}]$.
\par
  $(3)$ $\pi\neq (4^2,3^4)$,  $(4,3^4,2)$, $(5,4,3^5)$, $(5,3^5,2)$, $(4^7)$,
        $(4^3,3^4)$, $(4^2,3^4,2)$,\par  $(4,3^6)$, $(4,3^5,1)$,$(4,3^4,2^2)$, $(5,3^7)$,
        $(5,3^6,1)$, $(4^8)$, $(4^2,3^6)$, $(4^2,3^5,1)$,\par $(4,3^6,2)$,
        $(4,3^5,2,1)$, $(4,3^7,1)$, $(4,3^6,1^2)$,
        $(n-1,3^5,1^{n-6})$ and \par $(n-1,3^6,1^{n-7})$.
\par
{\bf Proof:} Assume that $\pi$ is potentially $K_5-2K_2$-graphic.
$(1)$ is obvious. According to Lemma 2.4, (2) holds.  Now it is easy
to check that $(4^2,3^4)$, $(4,3^4,2)$,\ \  $(5,4,3^5)$,\ \
$(5,3^5,2)$,\ \  $(4^7)$,\ \  $(4^3,3^4)$,\ \  $(4^2,3^4,2)$,\ \
$(4,3^6)$,\ \  $(4,3^5,1)$,\ \ $(4,3^4,2^2)$, $(5,3^7)$,
$(5,3^6,1)$, $(4^8)$, $(4^2,3^6)$, $(4^2,3^5,1)$, $(4,3^6,2)$,
$(4,3^5,2,1)$, $(4,3^7,1)$, $(4,3^6,1^2)$ are not potentially
$K_5-2K_2$-graphic and by Lemma 2.4, $\pi\neq(n-1,3^5,1^{n-6})$ and
$(n-1,3^6,1^{n-7})$. Hence, $(3)$ holds.
\par
  Now we prove the sufficient conditions. Suppose the graphic
sequence $\pi$ satisfies the conditions $(1)$-$(3)$. Our proof is by
induction on $n$. We first prove the base case where $n=5$. In this
case, $\pi$ is one of the following: $(4^5)$, $(4^3,3^2)$,
$(4,3^4)$. It is easy to check that all of these are potentially
$K_5-2K_2$-graphic. Now suppose that the sufficiency holds for
$n-1(n\geq6)$, we will show that $\pi$ is potentially
$K_5-2K_2$-graphic in terms of the following cases:
\par
\textbf{Case 1:} $d_n\geq4$. Clearly,
$\pi^\prime=(d_1^\prime,d_2^\prime,\cdots,d_n^\prime)$ satisfies
$(1)$ and (2). If $\pi^\prime$ also satisfies $(3)$, then by the
induction hypothesis, $\pi^\prime$ is potentially
$K_5-2K_2$-graphic, and hence so is $\pi$. If $\pi^\prime$ does not
satisfy $(3)$, since $\pi\neq(4^7)$ and $(4^8)$, then $\pi^\prime$
is just $(4^7)$ or $(4^8)$, and hence $\pi=(5^4,4^4)$ or
$(5^4,4^5)$. It is easy to check that these sequences are
potentially $K_5-2K_2$-graphic.
\par
\textbf{Case 2:} $d_n=3$. Consider
$\pi^\prime=(d_1^\prime,d_2^\prime,\cdots,d_{n-1}^\prime)$ where
$d_{n-3}^\prime\geq3$ and $d_{n-1}^\prime\geq2$. If $\pi^\prime$
satisfies $(1)$-$(3)$, then by the induction hypothesis,
$\pi^\prime$ is potentially $K_5-2K_2$-graphic, and hence so is
$\pi$.
\par
  If $\pi^\prime$ does not satisfy $(1)$, there are three subcases:
\par
\textbf{Subcase 1:} $d_1^\prime=d_5^\prime=3$. Then
$\pi^\prime=(3^k,2^{n-1-k})$ where $n-3\leq k\leq n-1$. Since
$\sigma(\pi^\prime)$ is even, $k$ must be even.  If $k=n-3$, then
$\pi=(4,3^{n-1})$ where $n$ is odd. Since $\pi\neq(4,3^6)$, we have
$n\geq9$. By Lemma 2.5, $\pi_1=(3^{n-5})$ is graphic. Let $G_1$ be a
realization of $\pi_1$, then $K_{1,2,2}\cup G_1$ is a realization of
$\pi=(4,3^{n-1})$. In other words, $\pi=(4,3^{n-1})$ is potentially
$K_5-2K_2$-graphic. If $k=n-2$, then $\pi=(4^2,3^{n-2})$ where $n$
is even. Since $\pi\neq(4^2,3^4)$ and $(4^2,3^6)$, we have $n\geq
10$. It is easy to see that $(4^2,3^8)$ and $(4^2,3^{10})$ are
potentially $K_5-2K_2$-graphic. Let $G_2$ be a realization of
$(4^2,3^8)$, which contains $K_5-2K_2$. If $n\geq14$, then
$\pi_2=(3^{n-10})$ is graphic by Lemma 2.5. Let $G_3$ be a
realization of $\pi_2$, then $G_2\cup G_3$ is a realization of
$\pi=(4^2,3^{n-2})$. In other words, $\pi=(4^2,3^{n-2})$ is
potentially $K_5-2K_2$-graphic. If $k=n-1$, then $\pi=(4^3,3^{n-3})$
where $n$ is odd. Since $\pi\neq(4^3,3^4)$, we have $n\geq9$.
Clearly, $K_5-e\cup G_1$ is a realization of $\pi=(4^3,3^{n-3})$.
Thus, $\pi=(4^3,3^{n-3})$ is potentially $K_5-2K_2$-graphic since
$K_5-2K_2\subseteq K_5-e$.
\par
\textbf{Subcase 2:} $d_1^\prime\geq 4$ and $d_5^\prime=2$. Since
$d_{n-3}^\prime \geq3$, we have $n=6$ or $n=7$. Then $\pi$ is
$(5^2,3^4)$, $(5,3^5)$ or $(6,3^6)$,  which is impossible by
condition (2) and (3).
\par
\textbf{Subcase 3:} $d_1^\prime=3$ and $d_5^\prime=2$. Then
$\pi=(4^2,3^4)$ or $(4,3^6)$, which is impossible by condition (3).
\par
  If $\pi^\prime$ does not satisfy $(2)$, then
$\pi^\prime=((n-2)^2,3^{n-3})$ or $((n-2)^2,3^{n-4},2)$. Hence,
$\pi=((n-1)^2,4,3^{n-3})$ or $((n-1)^2,3^{n-2})$. But
$\pi=((n-1)^2,3^{n-2})$ contradicts condition (2), thus
$\pi=((n-1)^2,4,3^{n-3})$.  Since $\pi_1^\prime=(n-2,3,2^{n-3})$ \ \
is \ \ potentially\ \  $C_4$-graphic\ \  by\ \  Theorem 2.3,\ \
thus\ \ $\pi=((n-1)^2,4,3^{n-3})$ is potentially $K_5-2K_2$-graphic.
\par
  If $\pi^\prime$ does not satisfy $(3)$, since $\pi\neq(5,4,3^5)$
and $(5,3^7)$, then $\pi^\prime$ is one of the following:
$(4^2,3^4)$, $(5,4,3^5)$, $(5,3^5,2)$, $(4^7)$, $(4^3,3^4)$,
$(4^2,3^4,2)$, $(4,3^6)$, $(5,3^7)$, $(4^8)$, $(4^2,3^6)$,
$(4,3^6,2)$, $(5,3^5)$, $(6,3^6)$. Hence, $\pi$ is one of the
following: $(5^2,4,3^4)$, $(5,4^3,3^3)$, $(4^5,3^2)$, $(6,5,4,3^5)$,
$(6,4^3,3^4)$, $(6,4,3^6)$, $(5^3,4^4,3)$, $(5^3,3^5)$,
$(5^2,4^2,3^4)$, $(5,4^4,3^3)$, $(4^6,3^2)$, $(5^2,3^6)$,
$(5,4^2,3^5)$, $(4^4,3^4)$, $(6,4^2,3^6)$, $(5^3,4^5,3)$,
$(5^2,4,3^6)$, $(5,4^3,3^5),$ $(4^5,3^4)$, $(5,4,3^7)$,
$(6,4^2,3^4)$, $(7,4^2,3^5)$. It is easy to check that all of these
are potentially $K_5-2K_2$-graphic.
\par
\textbf{Case 3:} $d_n=2$. Consider
$\pi^\prime=(d_1^\prime,d_2^\prime,\cdots,d_{n-1}^\prime)$ where
$d_4^\prime\geq3$ and $d_{n-1}^\prime\geq2$. If $\pi^\prime$
satisfies $(1)$-$(3)$, then by the induction hypothesis,
$\pi^\prime$ is potentially $K_5-2K_2$-graphic, and hence so is
$\pi$.
\par
  If $\pi^\prime$ does not satisfy $(1)$, there are three subcases:
\par
\textbf{Subcase 1:} $d_1^\prime=d_5^\prime=3$. Then $d_1=4$,
$d_3=d_4=d_5=3$ and $3\leq d_2\leq4$. If $d_2=4$, then
$\pi=(4^2,3^k,2^{n-2-k})$ where $k\geq3$ and $n-2-k\geq1$. Since
$\sigma(\pi)$ is even, $k$ must be even. We will show that $\pi$ is
potentially $K_5-2K_2$-graphic. It is enough to show
$\pi_1=(3^{k-3},2^{n-2-k},1)$ is graphic. If $n\geq8$, then $\pi_1$
is graphic by Lemma 2.5. If $n\leq7$, then $\pi=(4^2,3^4,2)$, which
is impossible by (3). If $d_2=3$, then $\pi=(4,3^k,2^{n-1-k})$ where
$k\geq6$, $n-1-k\geq1$ and $k$ is even. Since $\pi\neq(4,3^6,2)$, we
have $n\geq9$. We will show that $\pi$ is potentially
$K_5-2K_2$-graphic. It is enough to show $\pi_2=(3^{k-4},2^{n-1-k})$
is graphic. By Lemma 2.5, $\pi_2$ is graphic.
\par
\textbf{Subcase 2:} $d_1^\prime\geq 4$ and $d_5^\prime=2$. Then
$d_1\geq5$, $d_2=d_3=d_4=d_5=3$ and $d_6=\cdots=d_{n-1}=2$. Hence,
$\pi=(d_1,3^4,2^{n-5})$. Since $\sigma(\pi)$ is even, $d_1$ must be
even. We will show that $\pi$ is potentially $K_5-2K_2$-graphic. It
is enough to show $\pi_1=(d_1-4,2^{n-5})$ is graphic. It clearly
suffices to show $\pi_2=(2^{n-1-d_1},1^{d_1-4})$ is graphic. By
$\sigma(\pi_2)$ being even and Theorem 2.2, $\pi_2$ is graphic.
\par
\textbf{Subcase 3:} $d_1^\prime=3$ and $d_5^\prime=2$. Then
$\pi=(4,3^4,2^{n-5})$. Since $\pi\neq(4,3^4,2)$ and $(4,3^4,2^2)$,
we have $n\geq8$. Clearly, $K_{1,2,2}\cup C_{n-5}$ is a realization
of $\pi$. In other words, $\pi$ is potentially $K_5-2K_2$-graphic.
\par
  If $\pi^\prime$ does not satisfy $(2)$, i.e., $$ \pi^\prime= \left\{
    \begin{array}{ll}((n-2)^2,3^{n-3-2k}, 2^{2k}), \ \ \ \    \mbox{ $n$ is odd;}\\
    ((n-2)^2,3^{n-4-2k}, 2^{2k+1}),\ \
      \mbox{$n$ is even.} \end{array} \right. $$ If $n\geq7$, then $$ \pi= \left\{
    \begin{array}{ll}((n-1)^2,3^{n-3-2k}, 2^{2k+1}), \ \ \ \    \mbox{ $n$ is odd;}\\
    ((n-1)^2,3^{n-4-2k}, 2^{2k+2}),\ \ \ \ \ \
      \mbox{$n$ is even.} \end{array} \right.$$ which contradicts
condition (2). If $n=6$, then $\pi^\prime=(4^2,3^2,2)$ and hence
$\pi=(5^2,3^2,2^2)$ or $(4^4,2^2)$, which is impossible by (1).
\par
  If $\pi^\prime$ does not satisfy $(3)$, then $\pi^\prime$ is one of the following:
$(4^2,3^4)$,  $(4,3^4,2)$, $(5,4,3^5)$, $(5,3^5,2)$, $(4^7)$,
  $(4^3,3^4)$, $(4^2,3^4,2)$, $(4,3^6)$,  $(4,3^4,2^2)$, $(5,3^7)$,
   $(4^8)$, $(4^2,3^6)$, $(4,3^6,2)$, $(5,3^5)$, $(6,3^6)$. Since
   $\pi\neq(5,3^5,2)$, then $\pi$ is one of the
following: $(5^2,3^4,2)$, $(5,4^2,3^3,2)$, $(4^4,3^2,2)$,
$(5,4,3^3,2^2)$,\ \  $(4^3,3^2,2)$,\ \  $(6,5,3^5,2)$,\ \
 $(6,4^2,3^4,2)$,\ \
$(6,4,3^4,2^2)$,\ \  $(6,3^6,2)$,\ \  $(5^2,4^5,2)$,\ \
$(5^2,4,3^4,2)$, $(5,4^3,3^3,2)$,
 $(4^5,3^2,2)$, $(5^2,3^4,2^2)$,
$(5,4^2,3^3,2^2)$, $(4^4,3^2,2^2)$, $(5,4,3^5,2)$,\ \
$(4^3,3^4,2)$,\ \ $(5,4,3^3,2^3)$,\ \  $(5,3^5,2^2)$,\ \
$(4^3,3^2,2^3)$,\ \ $(6,4,3^6,2)$,\ \ $(5^2,4^6,2)$,\ \
$(5^2,3^6,2)$,\ \  $(5,4^2,3^5,2)$,\ \  $(4^4,3^4,2)$,\ \
$(5,4,3^5,2^2)$,\ \ $(5,3^7,2)$,\ \ $(4^3,3^4,2^2)$, $(6,4,3^4,2)$,
$(7,4,3^5,2)$. It is easy to check that all of these are potentially
$K_5-2K_2$-graphic.
\par
\textbf{Case 4:} $d_n=1$. Consider
$\pi^\prime=(d_1^\prime,d_2^\prime,\cdots,d_{n-1}^\prime)$ where
$d_5^\prime\geq3$. If $\pi^\prime$ satisfies $(1)$-$(3)$, then by
the induction hypothesis, $\pi^\prime$ is potentially
$K_5-2K_2$-graphic, and hence so is $\pi$.
\par
  If $\pi^\prime$ does not satisfy $(1)$, i.e., $d_1^\prime=3$, then
$d_1=4$ and $d_2=\cdots=d_5=3$. Hence, $\pi=(4,3^k,2^t,1^{n-1-k-t})$
where $k\geq4$ and $n-1-k-t\geq1$. Since $\sigma(\pi)$ is even,
$n-1-t$ must be even. We will show that $\pi$ is potentially
$K_5-2K_2$-graphic. It is enough to show
$\pi_1=(3^{k-4},2^t,1^{n-1-k-t})$ is graphic. Since
$\pi\neq(4,3^7,1)$ and $(4,3^6,1^2)$, we have $\pi_1\neq(3^3,1)$ and
$(3^2,1^2)$. If $n\geq9$, then $\pi_1$ is graphic by Lemma 2.5. If
$n\leq8$, since $\pi\neq(4,3^5,1)$ and $(4,3^5,2,1)$, then
$\pi=(4,3^4,1^2)$ or $(4,3^4,2,1^2)$. It is easy to see that $\pi$
is potentially $K_5-2K_2$-graphic.
\par
  If $\pi^\prime$ does not satisfy $(2)$, i.e., $$ \pi^\prime= \left\{
    \begin{array}{ll}(n-1-i,n-1-j,3^{(n-1)-i-j-2k}, 2^{2k},1^{i+j-2}), \\    \mbox{ $n-1-i-j$ is even;}\\
    (n-1-i,n-1-j,3^{(n-1)-i-j-2k-1}, 2^{2k+1},,1^{i+j-2}),\ \\
 \mbox{$n-1-i-j$ is odd.} \end{array} \right. $$ where $1\leq j\leq
     (n-1)-5$ and $0\leq k\leq [{{(n-1)-j-i-4}\over 2}]$. If $n-i>n-j+1$ or $n-i=n-j$, then $$ \pi= \left\{
    \begin{array}{ll}(n-i,n-(j+1),3^{n-i-(j+1)-2k}, 2^{2k},1^{i+(j+1)-2}),\\   \mbox{ $n-i-(j+1)$ is even;}\\
    (n-i,n-(j+1),3^{n-i-(j+1)-2k-1}, 2^{2k+1},1^{i+(j+1)-2}),
      \\\mbox{$n-i-(j+1)$ is odd.} \end{array} \right.$$ which contradicts
condition (2). If $n-i=n-j+1$, i.e., $$ \pi^\prime= \left\{
    \begin{array}{ll}(n-1-i,n-2-i,3^{n-2i-2k-2}, 2^{2k},1^{2i-1}),\\ \mbox{ $n$ is even;}\\
    (n-1-i,n-2-i,3^{n-2i-2k-3}, 2^{2k+1},1^{2i-1}),\ \
\\\mbox{$n$ is odd.} \end{array} \right. $$  Then $$ \pi= \left\{
    \begin{array}{ll}(n-i,n-i-2,3^{n-2i-2k-2}, 2^{2k},1^{2i}), \\    \mbox{ $n$ is even;}\\
    (n-i,n-i-2,3^{n-2i-2k-3}, 2^{2k+1},1^{2i}),\
      \\\mbox{$n$ is odd.}   \end{array} \right.$$ or $$ \pi= \left\{
    \begin{array}{ll}((n-1-i)^2,3^{n-2i-2k-2}, 2^{2k},1^{2i}), \\    \mbox{ $n$ is even;}\\
    ((n-1-i)^2,3^{n-2i-2k-3}, 2^{2k+1},1^{2i}),\\
      \mbox{$n$ is odd.}   \end{array} \right.$$ which contradicts
      condition (2).
\par
  If $\pi^\prime$ does not satisfy $(3)$, since $\pi\neq(5,3^6,1)$, $(4^2,3^5,1)$, $(n-1,3^5,1^{n-6})$
and $(n-1,3^6,1^{n-7})$, then $\pi^\prime$ is one of the following:
$(4^2,3^4)$,  $(4,3^4,2)$, $(5,4,3^5)$, $(5,3^5,2)$, $(4^7)$,
  $(4^3,3^4)$, $(4^2,3^4,2)$, $(4,3^5,1)$, $(4,3^4,2^2)$, $(5,3^7)$,
  $(5,3^6,1)$, $(4^8)$, $(4^2,3^6)$, $(4^2,3^5,1)$, $(4,3^6,2)$,
  $(4,3^5,2,1)$, $(4,3^7,1)$, $(4,3^6,1^2)$.
 Hence, $\pi$ is one of the following:
 \par
 $(5,4,3^4,1)$, $(4^3,3^3,1)$,
 $(5,3^4,2,1)$, $(4^2,3^3,2,1)$,
 $(6,4,3^5,1)$, $(5^2,3^5,1)$,
 \par
 $(6,3^5,2,1)$, $(5,4^6,1)$,
 $(5,4^2,3^4,1)$, $(4^4,3^3,1)$, $(5,4,3^4,2,1)$, $(4^3,3^3,2,1)$,
 \par
 $(5,3^5,1^2)$, $(4^2,3^4,1^2)$, $(5,3^4,2^2,1)$, $(4^2,3^3,2^2,1)$,
 $(6,3^7,1)$, $(6,3^6,1^2)$,\par
  $(5,4^7,1)$, $(5,4,3^6,1)$,
 $(4^3,3^5,1)$, $(5,4,3^5,1^2)$, $(4^3,3^4,1^2)$, $(5,3^6,2,1)$,\par
 $(4^2,3^5,2,1)$, $(5,3^5,2,1^2)$, $(4^2,3^4,2,1^2)$, $(5,3^7,1^2)$,
 $(4^2,3^6,1^2)$, $(5,3^6,1^3)$,\par
  $(4^2,3^5,1^3)$. It is
easy to check that all of these are potentially $K_5-2K_2$-graphic.
\par
\section{  Application }

\par
Using Theorem 3.1 and Theorem 3.3, we give simple proofs of the
following theorems due to  Lai:
\par
  \textbf{Theorem 4.1 }  (Lai [14])  For $n\geq5$,
$\sigma(K_5-P_3,n)=4n-4$.
\par
\textbf{Proof:} First we claim that for $n\geq5,
\sigma(K_5-P_3,n)\geq4n-4$. It is enough to show that there exists
$\pi_1$ with $\sigma(\pi_1)=4n-6$, such that $\pi_1$ is not
potentially $K_5-P_3$-graphic. Take $\pi_1=((n-1)^2,2^{n-2})$, then
$\sigma(\pi_1)=4n-6$, and it is easy to see that $\pi_1$ is not
potentially $K_5-P_3$-graphic by Theorem 3.1.
\par
  Now we show that if $\pi$ is an $n$-term $(n\geq5)$ graphical
sequence with $\sigma(\pi)\geq4n-4$, then there exists a realization
of $\pi$ containing $K_5-P_3$. Hence, it suffices to show that $\pi$
is  potentially $K_5-P_3$-graphic.

\par
  If $d_5=1$, then $\sigma(\pi)=d_1+d_2+d_3+d_4+(n-4)$ and
$d_1+d_2+d_3+d_4\leq12+(n-4)=n+8$. Therefore,
$\sigma(\pi)\leq2n+4<4n-4$, a contradiction. Thus, $d_5\geq2$.
\par
  If $d_3\leq2$, then $\sigma(\pi)\leq d_1+d_2+2(n-2)\leq2(n-1)+2(n-2)=4n-6<4n-4$, a contradiction.
Thus, $d_3\geq3$.

\par
  If $d_1\leq3$, then $\sigma(\pi)\leq3n<4n-4$, a contradiction.
Thus, $d_1\geq4$.

\par
  Since $\sigma(\pi)\geq4n-4$, then $\pi$ is not one of the
following: $(4,3^2,2^3)$, $(4,3^2,2^4)$, $(4,3^6)$. Thus, $\pi$
satisfies the conditions (1) and (2) in Theorem 3.1. Therefore,
$\pi$ is potentially $K_5-P_3$-graphic.
\par
\vspace{0.5cm}
  \textbf{Theorem 4.2 }  (Lai [13])  For $n\geq5$, $\sigma(K_5-C_4,n)=4n-4$.
  \par
 \textbf{Proof:} Obviously, for $n\geq5$, $\sigma(K_5-C_4,n)\leq
\sigma(K_5-P_3,n)=4n-4$. Now we claim $\sigma(K_5-C_4,n) \geq 4n-4$
for $n \geq5$. We would like to show there exists $\pi_1$ with
$\sigma(\pi_1)=4n-6$, such that $\pi_1$ is not potentially
$K_5-C_4$-graphic. Let $\pi_1=((n-1)^2,2^{n-2})$. It is easy to see
that $\sigma(\pi_1)=4n-6$ and the only realization of $\pi_1$ does
not contain $K_5-C_4$. Thus, $\sigma(K_5-C_4,n)=4n-4$.
\par
\vspace{0.5cm}
  \textbf{Theorem 4.3 }  (Lai [10], Luo[21])  $\sigma(C_5,n)=4n-4$ for $n\geq5$.
 \par
 \textbf{Proof:} Obviously, for $n\geq5$, $\sigma(K_5-C_5,n)\leq
\sigma(K_5-P_3,n)=4n-4$$(K_5-C_5=C_5)$. Now we claim $\sigma(C_5,n)
\geq 4n-4$ for $n \geq5$. We would like to show there exists $\pi_1$
with $\sigma(\pi_1)=4n-6$, such that $\pi_1$ is not potentially
$C_5$-graphic. Let $\pi_1=((n-1)^2,2^{n-2})$. It is easy to see that
$\sigma(\pi_1)=4n-6$ and the only realization of $\pi_1$ does not
contain $C_5$. Thus, $\sigma(C_5,n)=4n-4$.
\par
\vspace{0.5cm}
  \textbf{Theorem 4.4 }  (Lai [15])  For $n=5$ and $n\geq7$,
    $$\sigma(K_{3,1,1} ,n)=4n-2.$$
For $n=6$, if $\pi$ is a 6-term graphical sequence with $\sigma(\pi)
\geq 22$, then either there is a realization of $\pi$ containing
$K_{3,1,1}$ or $\pi=(4^{6})$. (Thus $\sigma(K_{3,1,1} ,6)=26$.)
\par
\textbf{Proof:} First we claim that for $n\geq5,
\sigma(K_5-K_3,n)\geq4n-2(K_{3,1,1}=K_5-K_3)$. It is enough to show
that there exists $\pi_1$ with $\sigma(\pi_1)=4n-4$, such that
$\pi_1$ is not potentially $K_5-K_3$-graphic. Take
$\pi_1=(n-1,3^{n-1})$, then $\sigma(\pi_1)=4n-4$, and it is easy to
see that $\pi_1$ is not potentially $K_5-K_3$-graphic by Theorem
3.3.
\par
  Now we show that if $\pi$ is an $n$-term $(n\geq5)$ graphical
sequence with $\sigma(\pi)\geq4n-2$, then there exists a realization
of $\pi$ containing $K_5-K_3$(unless $\pi=(4^6)$). Hence, it
suffices to show that $\pi$ is  potentially $K_5-K_3$-graphic.

\par
  If $d_5=1$, then $\sigma(\pi)=d_1+d_2+d_3+d_4+(n-4)$ and
$d_1+d_2+d_3+d_4\leq12+(n-4)=n+8$. Therefore,
$\sigma(\pi)\leq2n+4<4n-2$, a contradiction. Thus, $d_5\geq2$.
\par
  If $d_2\leq3$, then $\sigma(\pi)\leq d_1+3(n-1)\leq n-1+3(n-1)=4n-4<4n-2$, a contradiction.
Thus, $d_2\geq4$.
\par
  Since $\sigma(\pi)\geq4n-2$, then $\pi\neq(4^2,2^5)$. Hence, for $n=5$ and $n\geq7$, $\pi$ satisfies the conditions (1) and (2) in Theorem 3.3.
Therefore, $\pi$ is potentially $K_5-K_3$-graphic. For $n=6$, since
$\sigma(\pi)\geq4\times6-2=22$, then $\pi$ is not one of the
following: $(4^2,2^4)$, $(4^3,2^3)$. Thus, by Theorem 3.3, either
there is a realization of $\pi$ containing $K_{3,1,1}$ or
$\pi=(4^{6})$.
\par

        \section*{Acknowledgments}
  The authors are grateful to the referee for his valuable comments and suggestions.
\par

\end{document}